\documentclass[11pt,onecolumn]{article}
\usepackage[left=3cm,top=3cm,right=3cm,bottom=3cm]{geometry}
\usepackage{amsmath, verbatim, amssymb, color, amsthm,amsopn,mathrsfs,amsfonts}
\usepackage{graphicx,subfigure,grffile,times,caption,array,epstopdf}
\graphicspath{{./Figures/}}	
\newtheorem{theorem}{Theorem}
\newtheorem{ass}{Assumption}

\title{\vspace{-20pt}\Large{On the Existence of and Relationship between\\ Canards and Torus Canards in Forced Slow/Fast Systems}}
\author{
Han Wang\footnotemark[2],
Theodore Vo\footnotemark[2],
and Tasso J. Kaper\footnotemark[2] \footnotemark[3]
}

\begin{document}
\maketitle

\renewcommand{\thefootnote}{\fnsymbol{footnote}}
\footnotetext[2]{Department of Mathematics and Statistics, Boston University, 111 Cummington Mall, Boston, MA 02215, USA}  
\footnotetext[3]{\tt{tasso@math.bu.edu}}
\renewcommand{\thefootnote}{\arabic{footnote}}

\begin{abstract}
Canards are special solutions of slow/fast systems which are ubiquitous in neuroscience and electrical engineering. 
Two distinct classes of canard solutions have been identified and carefully studied: folded singularity canards and torus canards. 
Recently, an explicit and analytic relationship between these seemingly unrelated families of solutions was established in the classical forced van der Pol equation \cite{burke2015}. 
In this article, we generalize the results of \cite{burke2015} to the broader class of time-periodically forced planar slow/fast systems, which includes the forced van der Pol and the forced FitzHugh-Nagumo equations. 
We analytically determine the parameter values in this class of systems for which the two types of canard solutions exist, and show that the branches of primary canards of folded singularities continue into those of the torus canards as the forcing frequency is increased. 
We illustrate our results in the paradigm problem of the forced FitzHugh-Nagumo system.

	
\vspace{10pt} \noindent \textbf{Keywords:} folded singularity, relaxation oscillators, canards, torus canards, geometric desingularization, Melnikov theory, forced FitzHugh-Nagumo, forced van der Pol

\end{abstract}

\section{Introduction} 	\label{sec:intro}
	
Canards are solutions of slow/fast systems of differential equations that alternately spend long times near attracting and repelling states of the fast subsystem. 
Two distinct classes of canard solutions, distinguished by the underlying state of the fast subsystem, have been identified.
When the states are fast subsystem equilibria, the canard is of folded singularity type \cite{szmolyan2001}. 
When the states are limit cycles of the fast subsystem, the solution is known as a torus canard \cite{benes2011,burke2012,kramer2008}. 
Both folded singularity canards and torus canards are ubiquitous in applications. They have been used to explain the hard transition from small amplitude oscillations to relaxation oscillations in chemical reactions \cite{milik1998,petrov1992}, the pseudo-plateau bursting activity of neuroendocrine cells \cite{teka2011}, the mixed-mode oscillatory dynamics in hormone secretion by hypothalamic neurons in female mammals \cite{krupa2012}, the transient firing behaviour in neural excitability \cite{mitry2013,wechselberger2014}, and the transitions between tonic spiking and bursting states in paradigm neural bursting models \cite{burke2012}, among others.

\begin{figure}[h!]
\centering
\includegraphics[width=5in]{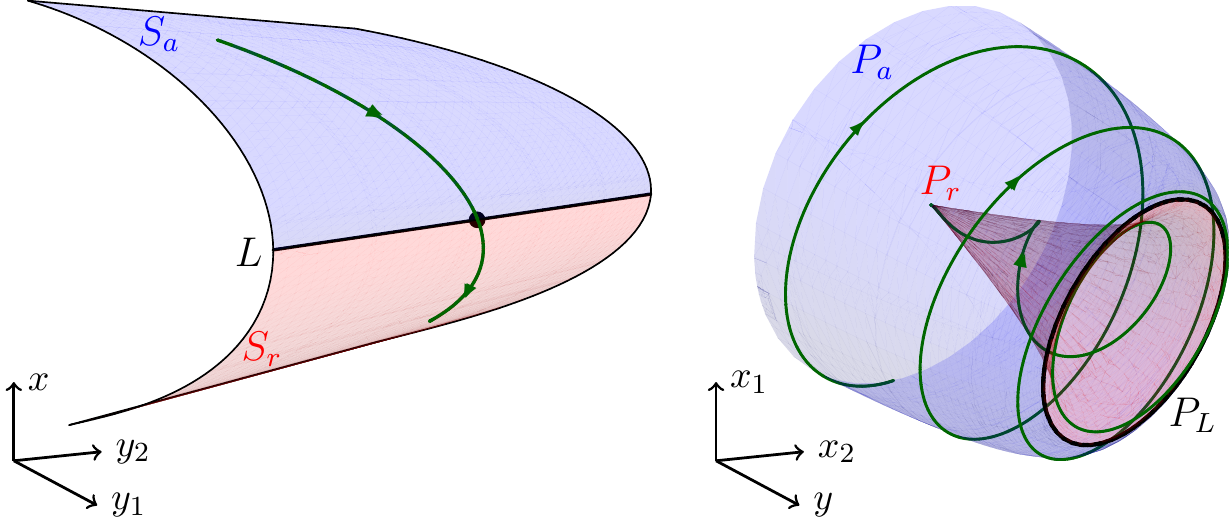}
\put(-362,145){(a)}
\put(-158,145){(b)}
\caption{Schematic of canard solutions in multi-scale systems with fast variables $x$ and slow variables $y$. (a) Folded singularity canards (green) connect attracting (blue) and repelling (red) manifolds, $S_a$ and $S_r$, of equilibria of the fast subsystem via a folded singularity (black marker). At least two slow variables are needed for these solutions to be generic. (b) Torus canard solutions (green) in slow/fast systems with at least two fast variables connect attracting (blue) and repelling (red) manifolds of limit cycles, $P_a$ and $P_r$, via a fold, $P_L$, of limit cycles.}
\label{fig:canards}
\end{figure}

Folded singularity canards alternate between attracting and repelling equilibria of the fast subsystem via a fold bifurcation of equilibria (Figure \ref{fig:canards}(a)). In slow/fast systems with only one slow variable, these canards are degenerate -- they require one-parameter families of $k$-fast/1-slow systems ($k \geq 1$) in order to be observed, and even then, they only occur on exponentially thin parameter sets \cite{diener1984,dumortier1996,krupa2001}. This exponential sensitivity is referred to as a canard explosion, and the
associated canard solutions are typically referred to as limit cycle
canards, or simply canard cycles \cite{BCDD1981, diener1984, dumortier1996}.
The addition of a second slow variable unfolds the canard phenomenon, making it generic and robust. 
That is, slow/fast systems with (at least) one fast variable and two slow variables can possess canard solutions on open parameter sets \cite{szmolyan2001, wechselberger2012}. 

An important aspect of canard theory is the behaviour of the slow manifolds of the system. 
Geometric singular perturbation theory (GSPT) \cite{fenichel1979,jones1995} guarantees that manifolds $S_a$ and $S_r$ of attracting and repelling fast subsystem equilibria persist as locally invariant attracting and repelling slow manifolds, $S_a^{\epsilon}$ and $S_r^{\epsilon}$, for sufficiently small perturbations, as measured by a small perturbation parameter $\epsilon$. GSPT breaks down in non-hyperbolic regions, such as neighbourhoods of a fast subsystem fold bifurcation. The extensions of $S_a^{\epsilon}$ and $S_r^{\epsilon}$ into neighbourhoods of a fold can result in unusual behaviour, such as local twisting of the slow manifolds. In particular, the slow manifolds can intersect in the neighbourhood of a special fold point, known as a folded singularity. These intersections between $S_a^{\epsilon}$ and $S_r^{\epsilon}$ are identified as maximal canard solutions \cite{szmolyan2001}. 

One of the remarkable results of canard theory \cite{szmolyan2001,wechselberger2005} is that the existence and bifurcations of maximal canards are encoded in the folded singularity itself. Four topologically distinct types of folded singularities, which are common in applications, are the folded node, folded saddle, and folded saddle-node of type I (FSN I) and type
II (FSN II). The folded node supports a finite and countable number of maximal canards (i.e., intersections between $S_a^{\epsilon}$ and $S_r^{\epsilon}$). The outermost and innermost intersections are known as the primary strong and weak maximal canards, respectively. The primary strong canard plays the role of a local phase space separatrix, dividing between those trajectories that rotate around the folded node, and those that do not. The primary weak canard is the axis of rotation for this local oscillatory behaviour. Additional maximal canards can exist between the primary canards, and these secondary canards further partition the slow manifolds based on the rotational properties of solutions. 

The folded saddle, on the other hand, has precisely one maximal canard solution and no rotational behaviour. This folded saddle canard typically divides the flow between solutions that transiently fire and those that remain quiescent \cite{wechselberger2014}. 
The FSN I is the codimension-one bifurcation in which the folded node and folded saddle coalesce and annihilate each other in a saddle-node bifurcation of folded singularities. The FSN I is known to possess maximal canards that inherit their properties from both the folded node and the folded saddle \cite{vo2015}. 
The FSN II is the codimension-one bifurcation in which a full system
equilibrium passes through and swaps stability with a folded singularity,
in a transcritical-type bifurcation. The FSN II has also been shown to
possess a family of maximal canards \cite{krupa2010}. We point out that in the planar case,
the FSN II is commonly referred to as a canard point, and the canards are
the limit cycle canards first discovered in \cite{BCDD1981, diener1984}.

	

The other main class of canard solutions, torus canards, was discovered in a model for the electrical activity of a cerebellar Purkinje cell \cite{kramer2008}. Torus canards are solutions that alternate between attracting and repelling manifolds of limit cycles of the fast subsystem via a saddle-node of limit cycles (Figure \ref{fig:canards}(b)). In analogy with folded singularity canards, maximal torus canards are defined as intersections between attracting and repelling invariant manifolds of limit cycles. 
Torus canards have been shown to occur ubiquitously in paradigm models from neuroscience and electrical engineering \cite{benes2011,burke2012,roberts2015}, where they typically mediate the transition between the tonic spiking and bursting regimes. The minimal dimensions required for a torus canard phenomenon are (at least) two fast variables and one slow variable. Note that when there is only one slow variable, the torus canards are degenerate, i.e., they only occur on exponentially thin parameter sets. Torus canards are generic in systems with two (or more) slow variables,
see \cite{vo2016}.


Recently, it was shown \cite{burke2015} that the classical forced van der Pol equation \cite{van1920, van1927, cartwright1950, cartwrigth1945, flaherty1978, levi1981, sekikawa2005} supports torus canard solutions. The method of geometric desingularization was used to analytically determine the parameter values --including forcing frequency, forcing amplitude, and the recovery threshold-- for which the torus canards exist. This parameter region coincides precisely with the transition region between the regimes of tonic spiking and bursting. Moreover, it was shown that as the forcing frequency decreased, the maximal torus canards would transition smoothly to primary maximal canards of FSN I type. This marked the first time that these two different types of canard solutions have been shown to be connected. The connection occurs smoothly as the forcing frequency is decreased and the states transition from being two-dimensional fast manifolds of persistent families of attracting and repelling limit cycles (for intermediate and high forcing frequencies) to being two-dimensional attracting and repelling slow manifolds (for low forcing frequencies).
	
In this article, we generalize the analytical results for the forced van der Pol equation \cite{burke2015} to the broader class of time-periodically forced planar slow/fast systems of the form
\begin{equation}	 \label{eqGen}
\begin{split}
x^\prime &= F\left( x,y,p,\epsilon \right),\\
y^\prime &= \epsilon \left( G \left( x,y,a,p,\epsilon \right) +b \cos\theta \right),\\
\theta^\prime &= \omega,
\end{split}
\end{equation}
where $(x,y) \in \mathbb{R}^2$, $a \in \mathbb{R}$ is a threshold parameter that is common to relaxation oscillators in electrical engineering and neuroscience, $p \in \mathbb{R}^k$ for $k \geq 1$ represents other system parameters, $0< \epsilon \ll 1$ is the ratio of slow and fast timescales, and $F$ and $G$ are $C^r$ with $r \geq 3$. 
We show, under fairly natural conditions on $F$ and $G$, that these general forced slow/fast systems have FSN I canards in the low frequency forcing regime (Theorem \ref{thmLowFreqF}), and torus canards in the intermediate forcing frequency regime (Theorem \ref{thmIntFreqF}). 
Moreover, we show that the maximal canards of the FSN I which exist for low frequency forcing continue into the maximal torus canards of the intermediate frequency forcing regime. 
Both theorems in this article are proven using the method of geometric desingularization (also known as the blow-up method) \cite{dumortier1993, dumortier1996}, and Melnikov theory, using techniques similar to those in \cite{burke2015}.
	
The article is organised as follows.
In Section \ref{sec:setup}, we outline the theoretical framework in which our analysis holds. 
In Section \ref{secLowFreqF}, we study folded singularity canards in the low frequency forcing regime (Theorem \ref{thmLowFreqF}). We follow in Section \ref{secIntFreqF} by analysing torus canards in the intermediate frequency regime (Theorem \ref{thmIntFreqF}). We then apply Theorems \ref{thmLowFreqF} and \ref{thmIntFreqF} to the forced FitzHugh-Nagumo equation \cite{fitzhugh1961,nagumo1962}, and thus present new analytical results for this 50-year old paradigm problem. 
We conclude in Section \ref{sec:lienard} where we show that, for small-amplitude forcing, our results are independent of whether the forcing enters via the slow direction or the fast direction. 

\section{Forced Planar Slow/Fast Systems Near Canard Points}  	\label{sec:setup}

Here, we state the assumptions of our problem in terms of the geometry of the unforced system
\begin{equation} \label{eqSlowSys}
\begin{aligned}
x' &= F\left( x,y,p,\epsilon \right),\\
y' &= \epsilon G \left( x,y,a,p,\epsilon \right).
\end{aligned}
\end{equation}
%

\begin{ass} 	\label{ass:recovery}
The system only depends linearly on the recovery threshold $a$.
\end{ass}

Taking the singular limit $\epsilon \rightarrow 0$ in (\ref{eqSlowSys}) gives the layer problem
\begin{equation} \label{eqLayer}
\begin{aligned}
x' &= F\left( x,y,p,0 \right),
\end{aligned}
\end{equation}
where the slow variable $y$ is a parameter. The critical manifold, $S$, is the set of equilibria of \eqref{eqLayer},
\begin{equation*}
S:= \left\{ \left( x,y\right) \in \mathbb{R}^2: F\left( x,y,p,0\right) = 0 \right\},
\end{equation*}
and is a key object in the geometric singular perturbations approach. In the unforced system \eqref{eqSlowSys}, the critical manifold is a 1D curve. Typically, as is often the case in applications, the interesting dynamics tend to manifest around non-hyperbolic regions of $S$, such as fold bifurcations of \eqref{eqLayer}. 

\begin{ass}	\label{ass:fold}
Let $B$ be an open subset of $\mathbb{R}^k$. For every $p \in B$, the critical manifold has a non-degenerate fold point at $\left(x_0(p), y_0(p)\right)$. That is,
\begin{equation}	\label{eq:fold}
\begin{split}
F\left( x_0, y_0, p, 0\right) &= 0, \qquad \frac{\partial F}{\partial x} \left(x_0, y_0, p, 0\right) = 0, \\
\frac{\partial^2 F}{\partial x^2} \left(x_0, y_0, p, 0\right) &\neq 0, \qquad \frac{\partial F}{\partial y} \left(x_0, y_0, p, 0\right) \neq 0.
\end{split}
\end{equation}
\end{ass}


\begin{ass}	\label{ass:canard}
System \eqref{eqSlowSys} possesses a non-degenerate canard point at $\mathscr{C} = \left( x_0, y_0, a_0, p, 0\right)$. That is, conditions \eqref{eq:fold} hold together with
\begin{equation}	\label{eq:canard} 
G\left( x_0, y_0, a_0, p, 0\right) = 0, \qquad
\frac{\partial G}{\partial x} \left(x_0, y_0, a_0, p, 0\right) \ne 0, \qquad 
\frac{\partial G}{\partial a} \left(x_0, y_0, a_0, p, 0\right) \ne 0.
\end{equation}
\end{ass}

The non-degeneracy condition $\left. \frac{\partial G}{\partial x} \right|_{\mathscr{C}} \neq 0$ ensures transverse intersection of the critical manifold $S$ and the slow nullcline $\{ G=0 \}$. The condition $\left. \frac{\partial G}{\partial a} \right|_{\mathscr{C}} \neq 0$ ensures that the intersection of $S$ and $\{ G=0 \}$ passes through the fold point with nonzero speed under variations in $a$ \cite{krupa2001}.

\begin{ass}	\label{ass:hopf}
The unforced system \eqref{eqSlowSys} possesses a Hopf bifurcation $\mathcal{O}(\epsilon)$-close to the canard point $\mathscr{C}$. Sufficient conditions for this are \eqref{eq:fold} and \eqref{eq:canard} together with
\[ \frac{\partial F}{\partial y} (x_0,y_0,p,0) \cdot \frac{\partial G}{\partial x} (x_0,y_0,a_0,p,0) < 0. \]
\end{ass}

We now consider planar slow/fast systems (within the framework of Assumptions \ref{ass:recovery}--\ref{ass:hopf}) subject to time-periodic forcing in the slow component, as given in \eqref{eqGen}.
%
%
The forcing frequency $\omega$ determines the type of slow/fast system and hence the nature of solutions that \eqref{eqGen} can support. For low frequencies ($\omega = \mathcal{O}(\epsilon)$), system \eqref{eqGen} is capable of generating folded singularity canards. For intermediate and high frequencies ($\omega = \mathcal{O}(\sqrt{\epsilon})$ and $\omega = \mathcal{O}(1)$), system \eqref{eqGen} potentially has torus canard solutions. 

There are many examples of these forced slow/fast systems in the literature. The most famous is the forced van der Pol equation \cite{van1920,van1927}
\begin{equation}	\label{eq:fvdp}
\begin{split}
x' &= y-\frac{x^3}{3} +x,\\
y' &= \epsilon \left( -x+a+b\cos\theta \right),\\
\theta ' &= \omega,
\end{split}
\end{equation}
which has the pair of non-degenerate canard points located at 
$ (x_0,y_0,a_0) = \pm \left( 1, - \frac{2}{3}, 1 \right)$.
%
Another prominent example is the forced FitzHugh-Nagumo system
\begin{equation} \label{eqFH-N}
\begin{split} 
x' &= x-\frac{x^3}{3} -y+I,\\
y' &= \epsilon \left( x+a-cy+b\cos\theta \right),\\
\theta ' &= \omega,
\end{split}
\end{equation}
which has non-degenerate canard points located at
$ \left( x_0,y_0,a_0 \right) =  \left( \pm 1, I \pm \frac{2}{3}, \left(I \pm \frac{2}{3} \right)c \mp 1 \right) $.
%

\section{Folded Singularity Canards For Low Forcing Frequencies} \label{secLowFreqF}

In the regime of low frequency forcing, we set $\omega = \epsilon \overline{\omega}$, where $\overline{\omega} = \mathcal{O} \left(1\right)$ with respect to $\epsilon$. System \eqref{eqGen} in that case is
\begin{equation} \label{eqGenLowFreqF}
\begin{split}
x' &= F\left( x,y,p,\epsilon \right),\\
y' &= \epsilon \left( G \left( x,y,a,p,\epsilon \right) +b\cos \theta \right),\\
\theta' &= \epsilon \overline{\omega},
\end{split}
\end{equation}
which potentially supports canard solutions. We first locate and classify the folded singularities of \eqref{eqGenLowFreqF} in order to study those canard solutions.

Taking the singular limit $\epsilon \to 0$ in (\ref{eqGenLowFreqF}) gives the layer problem, the dynamics of which are trivial on the critical manifold $S$. 
To describe the slow motions along $S$, we switch to the slow timescale ($\tau = \epsilon t$) and take the singular limit $\epsilon \to 0$ to obtain the reduced problem
\begin{equation} \label{eqReduced}
\begin{split}
0 &= F\left( x,y,p,0 \right),\\
\dot{y} &= G \left( x,y,a,p,0 \right) +b \cos\theta ,\\
\dot{\theta} &=  \overline{\omega},
\end{split}
\end{equation}
which describes the slow motions slaved to $S$. The reduced and layer problems are two different approximations of \eqref{eqGenLowFreqF}. The layer problem approximates the rapid motions away from the slow manifold, whilst the reduced problem approximates the slow motions restricted to the slow manifold. 
The idea of GSPT is that the dynamics of \eqref{eqGenLowFreqF} can be understood (for sufficiently small $\epsilon$) by suitably combining the dynamics of the reduced and layer problems.

\subsection{Folded Saddle-Nodes of Type I} 

The critical manifold of \eqref{eqGenLowFreqF} is a 2D surface, with attracting and repelling sheets, $S_a$ and $S_r$, separated by a curve of fold bifurcations (Assumption \ref{ass:fold}). The non-degeneracy conditions in Assumption \ref{ass:fold} guarantee that $S$ has a local graph representation, $y = y_S(x,p)$, in the neighbourhood of the fold curve. Using this graph representation, we can project the reduced flow onto the $(x,\theta)$-coordinate chart by differentiating the algebraic constraint in \eqref{eqReduced}, which gives  
%
%
\begin{equation} 	\label{RS}
\begin{split}
-\frac{\partial F}{\partial x}\dot{x} &=  \frac{\partial F}{\partial y} \left( G+b\cos\theta \right),\\
\dot{\theta} &= \overline{\omega},
\end{split}
\end{equation}
where all functions are evaluated along $S$, i.e., $y = y_S(x,p)$. We note that the reduced flow \eqref{RS} is singular along the fold curve. To remove this finite-time blow-up of solutions along the fold curve, we introduce the phase space-dependent time transformation, $dt =-\frac{\partial F}{\partial x} ds$, which gives the desingularized reduced system
\begin{equation} \label{eq:desing}
\begin{split}
\dot{x} &= \frac{\partial F}{\partial y} \left( G+b\cos\theta \right),\\
\dot{\theta} &= -\frac{\partial F}{\partial x} \overline{\omega},
\end{split}
\end{equation}
where the overdot now denotes derivatives with respect to $s$. The reduced and desingularized systems, \eqref{RS} and \eqref{eq:desing}, are topologically equivalent on $S_a$ (where $\frac{\partial F}{\partial x}<0$). On the repelling sheet $S_r$ (where $\frac{\partial F}{\partial x}>0$), the time transformation reverses the orientation of trajectories. In that case, the reduced flow is obtained from the desingularized flow by reversing the direction of orbits of \eqref{eq:desing}. 

The equilibria of the desingularized reduced system \eqref{eq:desing} are given by
\[ M := \left\{ (x,y,\theta) \in S : \frac{\partial F}{\partial x} =0, \quad G + b \cos \theta = 0 \right\}, \]
and are known as folded singularities. These are special points on the fold curve, where the $x$-equation of \eqref{RS} and \eqref{eq:desing} vanishes. That is, the folded singularities are points in the phase space where the reduced flow \eqref{RS} has a L'H\^{o}pital-type indeterminacy and solutions of \eqref{RS} can cross from $S_a$ to $S_r$ (or vice versa) with finite speed via the folded singularity. 

By Assumption \ref{ass:canard}, we have that system \eqref{eq:desing} possesses a pair folded singularities at $(x,y,a,p,\epsilon,\theta) = (x_0,y_0,a_0,p,0,\pm \frac{\pi}{2})$. Linear stability analysis of \eqref{eq:desing} shows that one of these is either a folded node or folded focus (depending on parameters), and the other is a folded saddle. Variations in the control parameter $a$ can alter the positions of these folded singularities. 
Taylor expanding $G$ about the canard point $\mathscr{C}$ and using the linear dependence on $a$ (Assumption \ref{ass:recovery}), we find that the pair of folded singularities merge to a FSN I singularity for the $a$ value such that 
\begin{align} \label{eq:FSNI} 
\left. \frac{\partial G}{\partial a} \right|_{\mathscr{C}} \left(a-a_0 \right) +b = 0. 
\end{align}
In order to state our results on the existence of primary canards of FSN I type in the context of \eqref{eqGenLowFreqF}, we introduce the convenient notation in Table \ref{tab:coefficients}, which will be needed in the Taylor series expansions of $F$ and $G$. 
We point out that Assumption \ref{ass:hopf} in this notation simply reads as $ c_2 c_3 < 0$.

\begin{table}[h!]
\centering
\renewcommand\arraystretch{2.5} 
\begin{tabular}{|c |c |c |c |c |c |c |c |} 
\hline
$c_1$ & $\displaystyle \frac{1}{2} \frac{\partial^2 F}{\partial x^2}$ &  
$c_2$ & $\displaystyle \frac{\partial F}{\partial y}$ & 
$c_3$ & $\displaystyle \frac{\partial G}{\partial x}$ & 
$c_4$ & $\displaystyle \frac{\partial G}{\partial a}$ \\
\hline
$c_5$ & $\displaystyle \frac{c_1}{c_2c_3} \frac{\partial F}{\partial \epsilon}$ &  
$c_6$ & $\displaystyle -\frac{1}{c_2c_3} \frac{\partial G}{\partial y}$ & 
$c_7$ & $\displaystyle -\frac{c_1}{c_2c_3^2} \frac{\partial G}{\partial \epsilon}$ & 
$c_8$ & $\displaystyle \frac{1}{c_1c_4} \frac{\partial^2 G}{\partial x \partial a}$  \\
\hline
$a_1$ & $\displaystyle -\frac{1}{c_1c_2} \frac{\partial^2 F}{\partial x \partial y}$ &  
$a_2$ & $\displaystyle -\frac{1}{6c_1^2} \frac{\partial^3 F}{\partial x^3}$ & 
$a_3$ & $\displaystyle -\frac{1}{c_2c_3} \frac{\partial^2 F}{\partial x \partial \epsilon}$ & 
$a_4$ & $\displaystyle -\frac{1}{2c_1c_3} \frac{\partial^2 G}{\partial x^2}$ \\
\hline
\end{tabular}
\caption{All functions are evaluated at the canard point at $\mathscr{C} = (x_0,y_0,a_0,p,0)$.}
\label{tab:coefficients}
\end{table}


\begin{theorem} \label{thmLowFreqF}
Let $b = \mathcal{O} (\sqrt{\epsilon})$ and $\overline{\omega} = \mathcal{O} (1)$. Then there exists an $\epsilon_0 >0$ such that for all $0< \epsilon< \epsilon_0$, there are two curves in the $(\omega,a)$ plane given by
\begin{equation} \label{eqthmLowFF}
a=a_0 -\frac{c_2 c_3^2}{c_1 c_4} \left( \frac{a_1}{8}-\frac{3}{8}a_2 + \frac{a_3}{2} + \frac{a_4}{4}-\frac{a_1c_5}{2}+c_5c_6 +\frac{c_6}{4} -c_7\right)\epsilon \pm \frac{b}{c_4} \exp \left(\frac{\epsilon \overline{\omega}^2}{2c_2c_3} \right) +\mathcal{O}(\epsilon^{3/2}),
\end{equation}
along which the system (\ref{eqGenLowFreqF}) has folds of primary maximal canards. Moreover, system (\ref{eqGenLowFreqF}) has two primary canards for every value of $a$ in the interval between these curves. There are no primary canards for values of $a$ outside the closure of these intervals.
\end{theorem}

\subsection{Preliminary Transformations and Blow-Up} 

In order to carry out the local analysis near the canard point, we translate the canard point to the origin via the coordinate transformation
\begin{equation} \label{eqCC1}
\begin{split}
& u = -c_1 \left( x-x_0 \right), \quad v = -c_1c_2 \left( y-y_0 \right), \quad \widetilde{a} = \frac{c_1c_4}{c_3} \left( a -a_0 \right), 
\end{split}
\end{equation}
and introduce the scaled parameters $(\widetilde{\epsilon}, \widetilde{b})$ according to
\begin{equation} \label{eqCC2}
\begin{split}
\widetilde{\epsilon} = -c_2c_3 \epsilon, \quad \widetilde{b}= \frac{c_1}{c_3}b.
\end{split}
\end{equation}
Taylor expanding $F$ and $G$ about the canard point $\mathscr{C}$, system \eqref{eqGenLowFreqF} becomes
\begin{equation} \label{eqOLowFreqF}
\begin{split}
u' &= vh_1-u^2h_2+\widetilde{\epsilon} h_3,\\
v' &= \widetilde{\epsilon} \left( -uh_4+vh_6+\widetilde{\epsilon} h_7+\widetilde{b}\left( \cos \theta-1\right)+\widetilde{b}+\widetilde{a}h_5 \right),\\
\theta' &= -\frac{\widetilde{\epsilon}}{c_2c_3} \overline{\omega},\\
\widetilde{\epsilon}' &= 0,
\end{split}
\end{equation}
where we have appended the trivial $\widetilde{\epsilon}$-equation to the system, and 
\begin{align*}
&h_j \left( u,v,\widetilde{a},\widetilde{\epsilon}\right) = 1+a_j u + \mathcal{O} \left( u^2,v,\widetilde{\epsilon}\right), \quad j = 1,2,\\
&h_3 \left( u,v,\widetilde{a},\widetilde{\epsilon}\right) = c_5+a_3 u + \mathcal{O} \left( u^2,v, \widetilde{\epsilon}\right),\\
&h_4 \left( u,v,\widetilde{a},\widetilde{\epsilon}\right) = 1+a_4u+c_8 \widetilde{a} + \mathcal{O} \left( u^2,v, \widetilde{\epsilon}\right)\\
&h_5 \left( u,v,\widetilde{a},\widetilde{\epsilon}\right) = 1+ \mathcal{O} \left( u^2,v,\widetilde{\epsilon}\right),\\
&h_k \left( u,v,\widetilde{a},\widetilde{\epsilon}\right) = c_k+ \mathcal{O} \left( u,v,\widetilde{a}, \widetilde{\epsilon}\right), \quad k=6,7.
\end{align*}

The origin is a nilpotent fixed point of \eqref{eqOLowFreqF}. In fact, there is a FSN I singularity at the origin whenever the parameters are such that
\[ \widetilde{a}+ \widetilde{b} = 0, \]
which is equivalent to the condition in \eqref{eq:FSNI}. 
Hence, we use the geometric desingularization method \cite{dumortier1993}, also known as the blow-up method, to analyze the dynamics around this point. The blow-up transformation is a map from the new variables $\left( \overline{u}, \overline{v}, \overline{\theta}, \overline{\epsilon}, \overline{r} \right) \in \mathbb{S}^3 \times \left[ -\mu, \mu \right]$, for sufficiently small $\mu$, to $\left( u,v, \theta, \widetilde{\epsilon} \right) \in \mathbb{R}^4$. This transformation enlarges the degenerate fixed point at the origin into a topological sphere, and hence effectively serves as a magnifying lens through which one can better analyze the system dynamics in the neighborhood of the degenerate point. 

For system (\ref{eqOLowFreqF}), the appropriate blow-up transformation is given by
\begin{equation}
u = \overline{r}^2 \overline{u}, \quad v = \overline{r}^4 \overline{v}, \quad \theta = \overline{r} \overline{\theta}, \quad \widetilde{\epsilon} = \overline{r}^4 \overline{\epsilon},
\end{equation}
which is commonly used in the blow-up analysis of FSN I singularities \cite{burke2015,vo2015}. 
In addition, it is useful to rescale the parameters as
\begin{equation} \label{eq:lowparamscalings}
\widetilde{a} = \sqrt{\widetilde{\epsilon}} \alpha, \quad \widetilde{b} = \sqrt{\widetilde{\epsilon}} \beta,  \quad \widetilde{a}+\widetilde{b} = \eta = \widetilde{\epsilon} \gamma,
\end{equation}
where $\alpha, \beta$, and $\gamma$ are $\mathcal{O} \left( 1 \right)$ with respect to $\widetilde{\epsilon}$.

As is common in these types of problems, one analyzes the flow on the topological hemisphere using an atlas of overlapping coordinate charts. 
The two most useful charts in the analysis of canard solutions \cite{dumortier1996,krupa2001,szmolyan2001} are the entry-exit chart $K_1 = \left\{\overline{v} = 1 \right\}$, and the rescaling chart $K_2 = \left\{ \overline{\epsilon} = 1\right\}$.

\subsection{The Entry-Exit Chart} 

In this chart, the coordinate change is given by
\[ u = r_1^2 u_1, \quad v = r_1^4, \quad \theta = r_1 \theta_1, \quad \widetilde{\epsilon} = r_1^4 \epsilon_1, \]
where the $1$-subscript indicates a variable in chart $K_1$. Transformation and desingularization (i.e., rescaling time by a factor of $r_1^2$) gives the blown-up system
\begin{align*}
\dot{u}_1 &= h_1-u_1^2 h_2+\epsilon_1 h_3 - \frac{1}{2} u_1 \epsilon_1 f,\\
\dot{r}_1 &= \frac{1}{4} r_1 \epsilon_1 f,\\
\dot{\theta}_1 &=-\frac{ r_1\epsilon_1}{c_2c_3} \overline{\omega} - \frac{1}{4} \theta_1 \epsilon_1 f,\\
\dot{\epsilon}_1 &= - \epsilon_1^2 f,
\end{align*}
where the overdot denotes derivatives with respect to the new time variable, and 
\[ f = -u_1h_4 +r_1^2 h_6 +r_1^2 \epsilon_1 h_7 +\sqrt{\epsilon_1} \beta \left( \cos \left(r_1 \theta_1\right)-1 \right) + \sqrt{\epsilon_1}\left( \beta+\alpha h_5\right). \]
The hyperplanes $\{r_1 = 0\}$ and $\{\epsilon_1 = 0 \}$ are invariant. The line
\[ l_u = \left\{ (u_1, r_1, \theta_1, \epsilon_1) = (u_1,0, 0, 0) \right\} \]
is invariant, and the system dynamics on this line are governed by $\dot{u}_1 = 1-u_1^2$ to leading order. There exist attracting and repelling fixed points $p_a = \left(1, 0, 0, 0\right)$ and $p_r = \left( -1, 0, 0, 0\right)$.
The fixed points, $p_a$ and $p_r$, have center manifolds, $N_{a,1}$ and $N_{r,1}$, respectively.

\subsection{The Rescaling Chart}

The canard solutions we seek enter the rescaling chart $K_2 = \left\{ \overline{\epsilon} = 1\right\}$, which is used to understand the dynamics along the top of the hemisphere. The coordinate change in this chart is
\[ u= r_2^2u_2, \quad v = r_2^4v_2, \quad \theta = r_2 \theta_2, \quad \widetilde{\epsilon} = r_2^4. \]
Thus, the blow-up transformation in $K_2$ is an $\epsilon$-dependent rescaling of the vector field \eqref{eqOLowFreqF}, i.e., it zooms in on an $\mathcal{O}(\widetilde{\epsilon}^{1/4})$-neighbourhood of the FSN I point.
Transformation and desingularization (by a factor of $r_2^2$) gives 
\begin{equation*}
\begin{split}
\dot{u}_2 &= v_2h_1-u_2^2h_2+h_3,\\
\dot{v}_2 &= -u_2h_4+r_2^2 v_2h_6 +r_2^2h_7+\beta \left( \cos \left( r_2\theta_2 \right)-1 \right) + \beta +\alpha h_5,\\
\dot{\theta}_2 &= -\frac{r_2}{c_2c_3} \overline{\omega},
\end{split}
\end{equation*}
where we again use the overdot to denote the derivative with respect to the new time variable.
We rewrite this third-order autonomous system as the second-order non-autonomous system
\begin{equation} \label{eq2ndNonAut}
\begin{split}
\dot{u}_2 = & v_2-u_2^2+c_5+ r_2^2 \left( a_1u_2v_2-a_2u_2^3+a_3u_2 \right) + \mathcal{O} \left( r_2^4 \right),\\
\dot{v}_2 = & -u_2+r_2^2 \left( -a_4u_2^2-c_8u_2\alpha+c_6v_2+c_7 +\gamma \right)\\ 
&+ \beta \left( \cos \left( -\frac{r_2^2\overline{\omega}}{c_2c_3}t_2\right) \cos \left(r_2 \theta_{2,0} \right)-1\right) -\beta \sin \left( -\frac{r_2^2\overline{\omega}}{c_2c_3}t_2\right) \sin \left(r_2 \theta_{2,0} \right)+ \mathcal{O} \left( r_2^4 \right),
\end{split}
\end{equation}
%
and note that
\begin{align*}
\cos \left( -\frac{r_2^2\overline{\omega}}{c_2c_3}t_2\right) \cos \left(r_2 \theta_{2,0} \right)-1 = \mathcal{O} \left( r_2^2 \right) & \text{ as } r_2 \rightarrow 0,\\
\sin \left( -\frac{r_2^2\overline{\omega}}{c_2c_3}t_2\right) \sin \left(r_2 \theta_{2,0} \right) = \mathcal{O} \left( r_2^3 \right) & \text{ as } r_2 \rightarrow 0.
\end{align*}

We analyze system \eqref{eq2ndNonAut} both in its unperturbed state (for $r_2=0$), and for small perturbations, i.e., in the regime where $r_2$ is small. The unperturbed version of system \eqref{eq2ndNonAut} is
\begin{equation} \label{eqUnperturbed}
	\begin{aligned}
	\dot{u}_2 &= v_2-u_2^2+c_5, \\
	\dot{v}_2 &= -u_2.
	\end{aligned}
\end{equation}
This is an integrable system with Hamiltonian function given by
\begin{equation}\label{eqHamiltonian}
	H \left( u_2,v_2 \right) = e^{-2v_2} \left( u_2^2-v_2-\frac{1}{2}-c_5\right).
\end{equation}
and non-canonical formulation
\begin{equation*}
	\begin{aligned}
	\dot{u}_2 &= \frac{1}{2} e^{2v_2} \frac{\partial H}{\partial v_2},\\
	\dot{v}_2 &= -\frac{1}{2} e^{2v_2} \frac{\partial H}{\partial u_2}.
	\end{aligned}
\end{equation*}

The special contour $\Gamma := \left\{ H=0 \right\}$ in chart $K_2$ is the heteroclinic on the upper hemisphere that connects the fixed points $p_r$ and $p_a$ (in chart $K_1$) on the equator of the hemisphere. The contour $\Gamma$ has the explicit time parametrization 
\begin{equation} \label{eqTimeParametrization}
\left(u_{2,\Gamma}, v_{2,\Gamma} \right) = \left( -\frac{1}{2} t_2, \frac{1}{4} t_2^2 -\frac{1}{2}-c_5 \right).
\end{equation}
This parabola also separates closed periodic orbits $\left( H<0 \right)$ from unbounded orbits $\left( H>0 \right)$, and corresponds to the strong canard of the FSN I. 

\subsection{Proof of Theorem \ref{thmLowFreqF}}

We now analyze \eqref{eq2ndNonAut} for small $r_2$ perturbations using Melnikov theory, which measures the splitting distance $D$ between the curves of solutions that are forward and backward asymptotic to the points $p_r$ and $p_a$. That is, the Melnikov method determines the parameter values for which the heteroclinic $\Gamma$ persists for $0<r_2 \ll 1$.
The Melnikov function $D$ has the asymptotic series expansion 
\[  D\left( r_2 \right) = d_{r_2^2} r_2^2 + d_{r_2^3} r_2^3 + \mathcal{O}\left( r_2^4 \right), \]
where the Melnikov integrals are given by
\begin{equation*}
\begin{aligned}
d_{r_2^2} &= \int_{-\infty}^{\infty} \nabla H|_{\Gamma} \cdot \left. \left( \begin{matrix}
a_1u_2v_2-a_2u_2^3+a_3u_2\\
-a_4u_2^2-c_8u_2\alpha+c_6v_2+c_7 +\gamma+\frac{\beta}{r_2^2} \left( \cos \left( -\frac{r_2^2\overline{\omega}}{c_2c_3}t_2\right) \cos \left(r_2 \theta_{2,0} \right)-1\right) 
\end{matrix} \right)\right|_{\Gamma} dt_2, \\
&= e^{1+2c_5} \sqrt{2\pi} \left( \frac{a_1}{8}-\frac{3}{8}a_2 + \frac{a_3}{2} + \frac{a_4}{4}-\frac{a_1c_5}{2} +c_5c_6 +\frac{c_6}{4}-c_7-\gamma -\frac{\beta}{r_2^2} \left( e^{-\frac{r_2^4 \overline{\omega}^2}{2c_2^2c_3^2}} \cos \left( r_2\theta_{2,0}\right) -1 \right) \right).\\
d_{r_2^3} &= \int_{-\infty}^{\infty} \nabla H|_{\Gamma} \cdot \left. \left( \begin{matrix}
0\\
- \frac{\beta}{r_2^2} \sin \left( -\frac{r_2^2\overline{\omega}}{c_2c_3}t_2\right) \sin \left(r_2 \theta_{2,0} \right)
\end{matrix} \right)\right|_{\Gamma} dt_2 = 0.
\end{aligned}
\end{equation*}
Substituting into the bifurcation equation $D \left( r_2 \right) = 0$, we find 
\begin{equation}
\widetilde{a} = \left( \frac{a_1}{8}-\frac{3}{8}a_2 + \frac{a_3}{2} + \frac{a_4}{4}-\frac{a_1c_5}{2} +c_5c_6 +\frac{c_6}{4}-c_7\right)\widetilde{\epsilon} -\widetilde{b} {\rm exp} \left(-\frac{\widetilde{\epsilon} \overline{\omega}^2}{2c_2^2c_3^2}\right) \cos \left(\theta_0\right) +\mathcal{O}(\widetilde{\epsilon}^{3/2}),
\end{equation}
where $\theta_0 = r_2 \theta_{2,0}$ is an arbitrary phase.
Reverting to the original parameters using \eqref{eqCC1}, \eqref{eqCC2}, and \eqref{eq:lowparamscalings}, we obtain
\begin{equation} \label{eq:fnc}
\begin{split}
a&=a_0 -\frac{c_2 c_3^2}{c_1 c_4} \left( \frac{a_1}{8}-\frac{3}{8}a_2 + \frac{a_3}{2} + \frac{a_4}{4}-\frac{a_1c_5}{2} +c_5c_6 +\frac{c_6}{4}-c_7\right)\epsilon -\frac{b}{c_4} \cos \theta_0 {\rm exp} \left(\frac{\epsilon \overline{\omega}^2}{2c_2c_3} \right)+\mathcal{O}(\epsilon^{3/2}).\\
\end{split}
\end{equation}
The envelope of these curves is exactly 
formula \eqref{eqthmLowFF}.
Therefore, we have shown that the heteroclinic $\Gamma$ persists along the curves \eqref{eqthmLowFF} in parameter space, and these correspond to maximal canard solutions of FSN I type for system \eqref{eqGenLowFreqF}. 
For each $a$ between the curves \eqref{eqthmLowFF}, there are primary maximal canards given by \eqref{eq:fnc}; and, for any value of $a$ outside these two curves there are no folded singularity canard solutions, since $\Gamma$ does not persist for these values of $a$. 
This concludes the proof of Theorem \ref{thmLowFreqF}.

\section{Torus Canards for Intermediate Forcing Frequencies} \label{secIntFreqF}

In this section, we analyze system \eqref{eqGen} in the regime of intermediate frequency forcing, where we set $\omega = \sqrt{\epsilon} \Omega$, with $\Omega = \mathcal{O} \left( 1 \right)$ with respect to $\epsilon$. System (\ref{eqGen}) becomes
\begin{equation} \label{eqGenIntFreqF}
\begin{aligned}
x' &= F\left( x,y,p,\epsilon \right),\\
y' &= \epsilon \left( G \left( x,y,a,p,\epsilon \right) +b \cos \theta \right),\\
\theta' &= \sqrt{\epsilon} \Omega.
\end{aligned}
\end{equation}
Recall by Assumption \ref{ass:hopf} that $c_2 c_3 <0$.

\begin{theorem} \label{thmIntFreqF}
	Let $b = \mathcal{O} (\epsilon)$ and $\Omega= \mathcal{O}(1)$. Then there exists an $\epsilon_0 >0$ such that for all $0< \epsilon < \epsilon_0$, there are two curves in the $(\omega, a)$ plane given by
\begin{equation} \label{eqthmIntFF}
a = a_0 -\frac{c_2 c_3^2}{c_1 c_4}  \left( \frac{a_1}{8}-\frac{3}{8}a_2 + \frac{a_3}{2} + \frac{a_4}{4}-\frac{a_1c_5}{2}+c_5c_6 +\frac{c_6}{4} -c_7 \right) \epsilon \pm \frac{b}{c_4}  \exp \left(\frac{\Omega^2}{2c_2c_3} \right)+\mathcal{O}(\epsilon^{3/2}),
\end{equation}
along which the system (\ref{eqGenIntFreqF}) has folds of maximal torus canards. Moreover, the system (\ref{eqGenIntFreqF}) has two torus canards for every value of $a$ in the interval between these curves. There are no torus canards for values of $a$ outside the closure of these intervals.
\end{theorem}

Formula \eqref{eqthmIntFF} of Theorem \ref{thmIntFreqF} is equivalent to formula \eqref{eqthmLowFF} of Theorem \ref{thmLowFreqF}. 
In fact, \eqref{eqthmLowFF} and \eqref{eqthmIntFF} are different representations (in different frequency regimes) of the same unifying formula
\begin{equation*} 
a = a_0 -\frac{c_2 c_3^2}{c_1 c_4}  \left( \frac{a_1}{8}-\frac{3}{8}a_2 + \frac{a_3}{2} + \frac{a_4}{4}-\frac{a_1c_5}{2}+c_5c_6 +\frac{c_6}{4} -c_7 \right) \epsilon \pm \frac{b}{c_4} \exp \left(\frac{\omega^2}{2c_2c_3 \epsilon} \right)+\mathcal{O}(\epsilon^{3/2}),
\end{equation*}
which confirms that the primary maximal canards of the FSN I in the low frequency regime continue smoothly into the maximal torus canards of the intermediate frequency regime.
However, the analysis needed to establish the validity of \eqref{eqthmIntFF} is different from that needed for \eqref{eqthmLowFF} due to the higher forcing frequency and different scalings of the parameters $a$ and $b$, and the structural difference in applying the blow-up technique.

Again, we translate and rescale the dependent variables so that the fold curve coincides with the $\theta$-axis. Here, we recall \eqref{eqCC1} and \eqref{eqCC2}, and the definitions of the functions $h_i$. In this intermediate frequency forcing regime, we rescale the parameters as
\[  \widetilde{a} = \widetilde{\epsilon}\,  \widetilde{\alpha}, \quad \widetilde{b} = \widetilde{\epsilon}\, \widetilde{\beta}, \]
where $\widetilde{\alpha}$, and $ \widetilde{\beta}$ are $\mathcal{O} \left( 1\right)$ with respect to $\widetilde{\epsilon}$. The governing equations become
\begin{equation} \label{eq:intermediate}
\begin{aligned}
u' &= v h_1 -u^2 h_2 + \widetilde{\epsilon} h_3,\\
v' &= \widetilde{\epsilon} \left( -u h_4 +v h_6 \right) + \widetilde{\epsilon}^2 \left(\widetilde{\alpha} h_5 + h_7 + \widetilde{\beta} \cos \theta \right),\\
\theta' &= \sqrt{-\frac{\widetilde{\epsilon}}{c_2c_3}}\Omega,\\
\widetilde{\epsilon}' &=0.
\end{aligned}
\end{equation}

Next, we use the geometric desingularization method to carry out the analysis near the circle of fold points. For all $\theta \in \mathbb{S}^1$, the transformation is a map from $\left( \overline{u}, \overline{v}, \overline{\epsilon}, \overline{r} \right) \in \mathbb{S}^2 \times \left[ -\mu, \mu \right] $ to $\left( u, v, \widetilde{\epsilon} \right) \in \mathbb{R}^3$. For system \eqref{eq:intermediate}, the coordinate change is given by
\[ u = \overline{r} \overline{u}, \quad v = \overline{r}^2 \overline{v}, \quad \widetilde{\epsilon} = \overline{r}^2 \overline{\epsilon}. \]
This is a cylindrical blow-up transformation which blows up the circle of fold points into a 2-torus (as opposed to the spherical blow-up in the previous section). Again, we focus on the entry-exit and rescaling charts. 

\subsection{The Entry-Exit Chart}

In the entry-exit chart $K_1= \left\{ \overline{v} = 1 \right\}$, the coordinate change is given by 
\begin{equation*}
u = r_1 u_1, \quad v = r_1^2, \quad \widetilde{\epsilon} = r_1^2 \epsilon_1.
\end{equation*}
In these blown-up coordinates, system \eqref{eq:intermediate} becomes
\begin{align*}
\dot{u}_1 &= h_1-u_1^2 h_2+\epsilon_1 h_3 - \frac{1}{2} u_1 \epsilon_1 f,\\
\dot{r}_1 &= \frac{1}{2} r_1 \epsilon_1 f,\\
\dot{\theta} &= \sqrt{\frac{\epsilon_1}{-c_2c_3}} \Omega,\\
\dot{\epsilon}_1 &= - \epsilon_1^2 f.
\end{align*}
where $f = -u_1h_4 +r_1 \left( \epsilon_1 \widetilde{\alpha} h_5 + \epsilon_1 \widetilde{\beta} \cos \theta+ h_6+\epsilon_1 h_7 \right)$, we have rescaled time by a factor of $r_1$, and the overdot denotes the derivative with respect to the new time variable. The hyperplanes $\left\{ r_1 = 0 \right\}$ and $\left\{ \epsilon_1 = 0 \right\}$ are invariant. The line
\[ l_u = \left\{ (u_1, r_1, \epsilon_1, \theta) = (u_1,0, 0, \theta) \right\} \]
is invariant, and the system dynamics on this line are governed by
\begin{align*}
\dot{u}_1 &= 1-u_1^2,\\
\dot{\theta} &= 0.
\end{align*}
There exist attracting and repelling fixed points $p_a = \left( 1, 0, 0, \theta \right)$ and $p_r = \left( -1, 0, 0, \theta \right)$. 
The fixed points have center manifolds $N_{a,1}$ and $N_{r,1},$ respectively.

\subsection{The Rescaling Chart}

In the rescaling chart $K_2 = \left\{ \overline{\epsilon} = 1\right\}$, the coordinate change is given by
\[ u=r_2u_2, \quad v = r_2^2v_2, \quad \widetilde{\epsilon} = r_2^2. \]
Transformation and desingularization (by a factor of $r_2$) gives the blown-up system
\begin{equation} \label{eq3rdInt}
\begin{split}
\dot{u}_2 &= v_2h_1-u_2^2h_2+h_3,\\
\dot{v}_2 &= -u_2h_4+r_2 \left( \widetilde{\alpha} h_5 + \widetilde{\beta} \cos \theta +v_2h_6 +h_7\right),\\
\dot{\theta} &= \frac{\Omega}{\sqrt{-c_2c_3}},\\
\dot{r}_2 &= 0,
\end{split}
\end{equation}
where the overdot denotes derivatives with respect to the new time. For the Melnikov analysis, it is convenient to rewrite (\ref{eq3rdInt}) as a second order non-autonomous system
\begin{equation} \label{eq2ndInt}
\begin{aligned}
\dot{u}_2 &= v_2-u_2^2+c_5+ r_2 \left( a_1u_2v_2-a_2u_2^3+a_3u_2 \right) + \mathcal{O} \left( r_2^2 \right),\\
\dot{v}_2 &= -u_2+r_2 \left( -a_4u_2^2+ \widetilde{\alpha} +c_6 v_2 +c_7+ \widetilde{\beta} \cos \left( \frac{\Omega}{\sqrt{-c_2c_3}} t_2 +\theta_0 \right) \right) + \mathcal{O} \left( r_2^2 \right).\\
\end{aligned}
\end{equation}
The unperturbed version of (\ref{eq2ndInt}) is exactly system (\ref{eqUnperturbed}), which is Hamiltonian, and it has the same heteroclinic orbit $\Gamma$, which is parametrized by \eqref{eqTimeParametrization}.

\subsection{Proof of Theorem \ref{thmIntFreqF}}

As in the proof of Theorem \ref{thmLowFreqF}, we use the Melnikov method to determine the parameter values for which the heteroclinic orbit $\Gamma$ persists for $0<r_2 \ll 1$.
The splitting distance in this case has the series expansion
\[ D\left( r_2 \right) = d_{r_2} r_2 + \mathcal{O}(r_2^2), \]
where the Melnikov integrals are given by
\begin{align*}
d_{r_2} &= \int_{-\infty}^{\infty} \nabla H|_{\Gamma} \cdot \left. \left( \begin{array}{c}
a_1u_2v_2-a_2u_2^3+a_3u_2\\
-a_4u_2^2+\widetilde{\alpha}+c_6v_2 +c_7+\widetilde{\beta} \cos \left( \frac{\Omega}{\sqrt{-c_2c_3}} t_2 +\theta_0 \right) 
\end{array} \right) \right|_{\Gamma} dt_2\\
&= e^{1+2c_5} \sqrt{2\pi} \left( \frac{a_1}{8}-\frac{3}{8}a_2 + \frac{a_3}{2} + \frac{a_4}{4}-\frac{a_1c_5}{2}+c_5c_6 +\frac{c_6}{4} -c_7-\widetilde{\alpha} -\widetilde{\beta} e^{\frac{\Omega^2}{2c_2c_3}} \cos\theta_0 \right).
\end{align*}
Substituting into the bifurcation equation $D\left( r_2 \right) =0$, we find that
\begin{equation}	\label{eq:tc}
\widetilde{a} = \left( \frac{a_1}{8}-\frac{3}{8}a_2 + \frac{a_3}{2} + \frac{a_4}{4}-\frac{a_1c_5}{2}+c_5c_6 +\frac{c_6}{4} -c_7 \right) \widetilde{\epsilon} - \widetilde{b} \cos \theta_0 \exp \left(\frac{\Omega^2}{2c_2c_3} \right) +\mathcal{O}(\widetilde{\epsilon}^{3/2}).
\end{equation}
Upon reverting to the original parameters, we obtain
\begin{align*}
a&=a_0 -\frac{c_2 c_3^2}{c_1 c_4} \left( \frac{a_1}{8}-\frac{3}{8}a_2 + \frac{a_3}{2} + \frac{a_4}{4}-\frac{a_1c_5}{2}+c_5c_6 +\frac{c_6}{4} -c_7 \right) \epsilon -\frac{b}{c_4} \cos \theta_0 \exp \left(\frac{\Omega^2}{2c_2c_3} \right)+\mathcal{O}(\epsilon^{3/2}).
\end{align*}
The envelope of these curves is precisely \eqref{eqthmIntFF}.
Therefore, we have shown that the heteroclinic $\Gamma$ persists along the curves \eqref{eqthmIntFF} in parameter space, and these correspond to maximal torus canard solutions of \eqref{eqGenIntFreqF}. 
Also, for each value of $a$ in between the two curves \eqref{eqthmIntFF}, there is a pair of torus canards given by \eqref{eq:tc}; and, for any value of $a$ outside these two curves there are no torus canard solutions, since $\Gamma$ does not persist for these values of $a$. This concludes the proof of Theorem \ref{thmIntFreqF}.

\section{Application to the forced FitzHugh-Nagumo equation} \label{secApp}

In this section, we apply Theorems \ref{thmLowFreqF} and \ref{thmIntFreqF} to the forced FitzHugh-Nagumo equation \eqref{eqFH-N}. 
Note that $p = (I,c)$ in this case. 
Here, we focus on the behaviour of the forced system near the non-degenerate canard point at  $(x_0,y_0,a_0) = (1,I+\frac{2}{3},-1+\left(\frac{2}{3}+I\right) c)$.
The coefficients for system \eqref{eqFH-N} about this canard point, as calculated
from the general formulas given in Table \ref{tab:coefficients}, are given in Table \ref{tab:FHN}. 
%
\begin{table}[h!]
\centering
\begin{tabular}{|r |c |r |c |r |c |r |c |}
\hline
$c_1$ & $-1$ & $c_2$ & $-1$ & $c_3$ & $1$ & $c_4$ & $1$ \\
\hline
$c_5$ & $0$ & $c_6$ & $-c$ & $c_7$ & $0$ & $c_8$ & $0$ \\
\hline
$a_1$ & $0$ & $a_2$ & $1/3$ & $a_3$ & $0$ & $a_4$ & $0$ \\
\hline
\end{tabular}
\caption{Coefficients for system \eqref{eqFH-N} about the canard point $(x_0,y_0,a_0) = (1,I+\frac{2}{3},-1+\left(\frac{2}{3}+I\right) c)$.}
\label{tab:FHN}
\end{table}
%

%
%

\begin{figure}[h!]
\centering
\subfigure[$I=0$]{\includegraphics[width=2.475in]{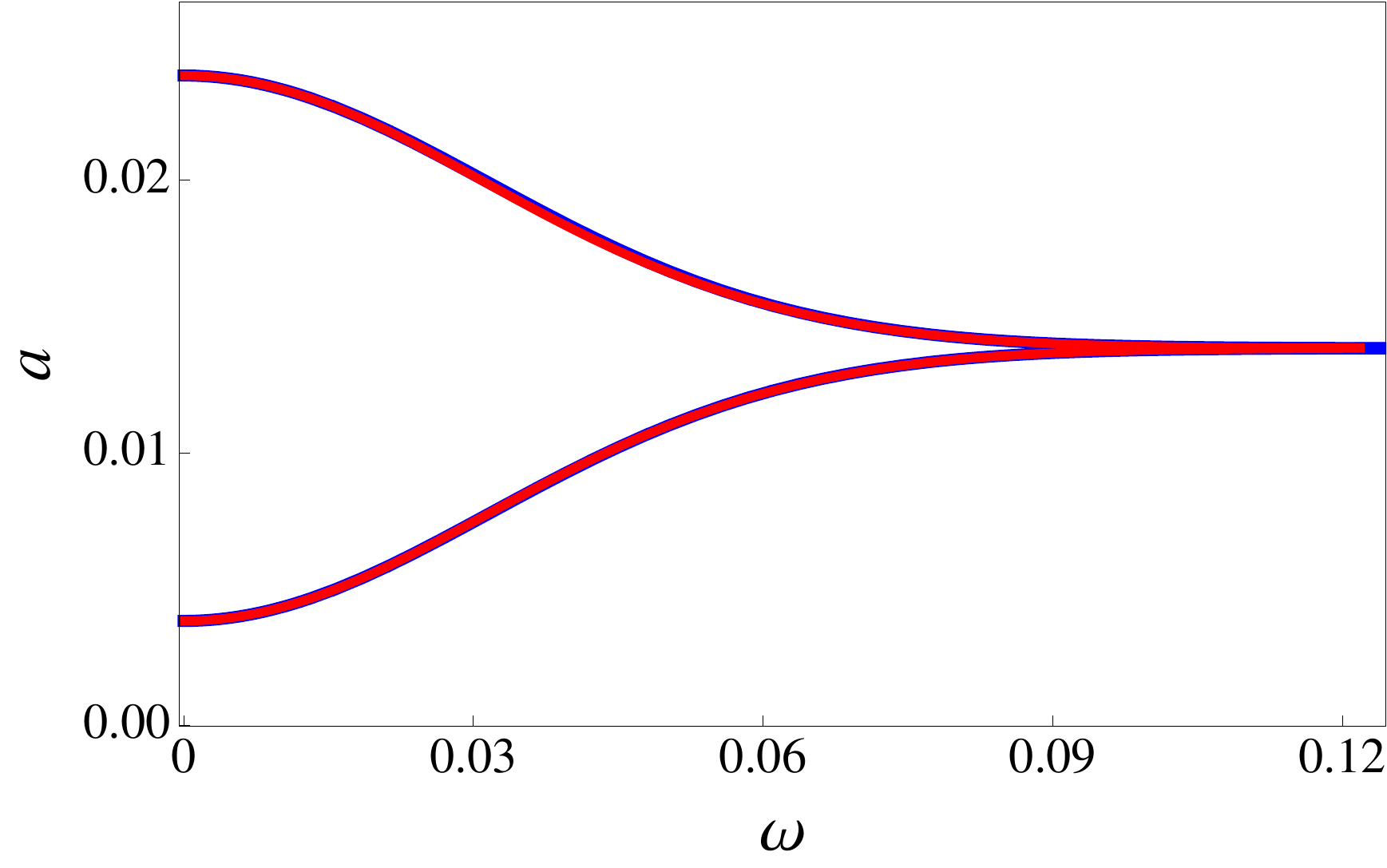}} \hspace{10pt}
\subfigure[$I=0.0015$]{\includegraphics[width=2.475in]{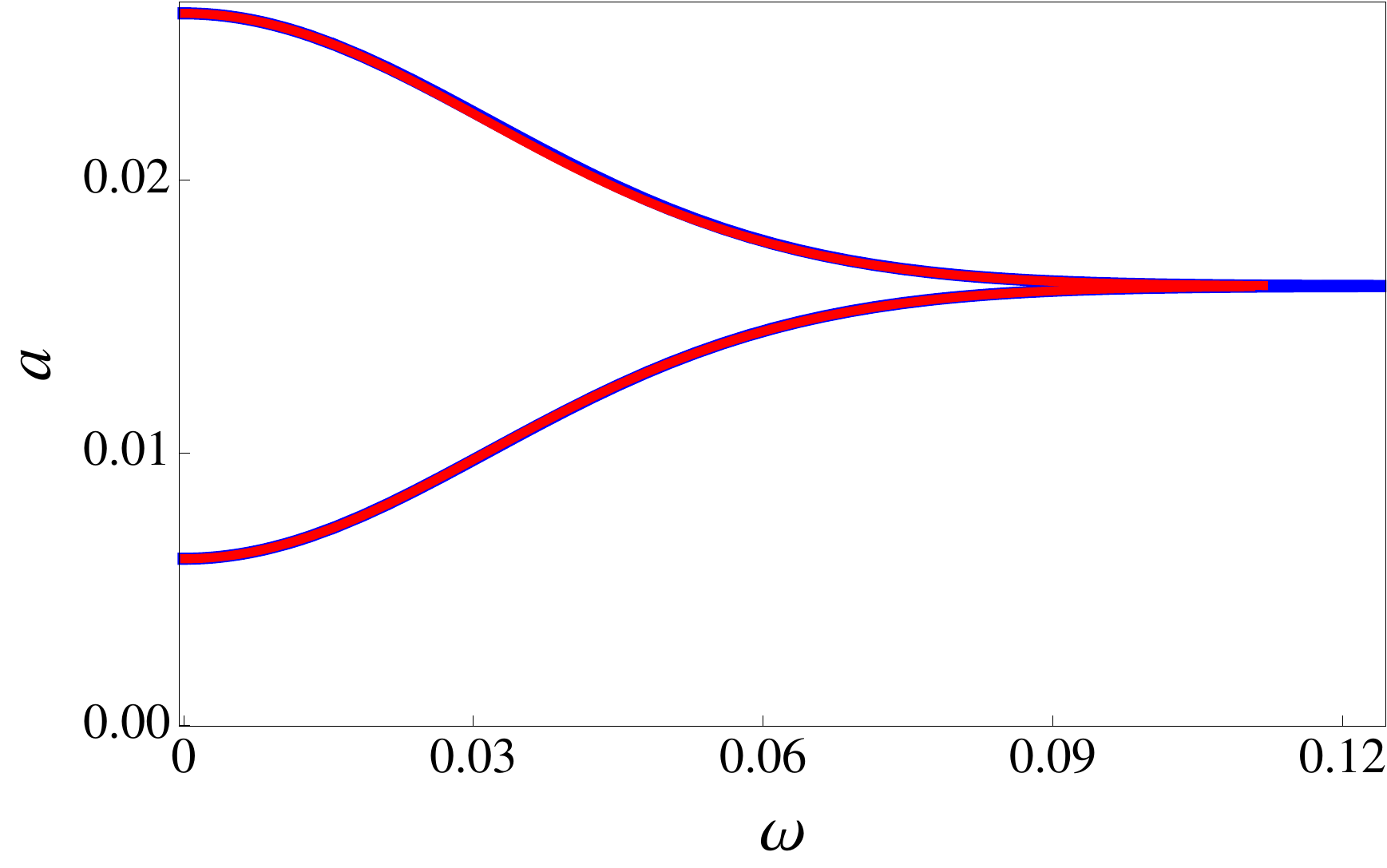}} \\
\subfigure[$c=1.515$]{\includegraphics[width=2.475in]{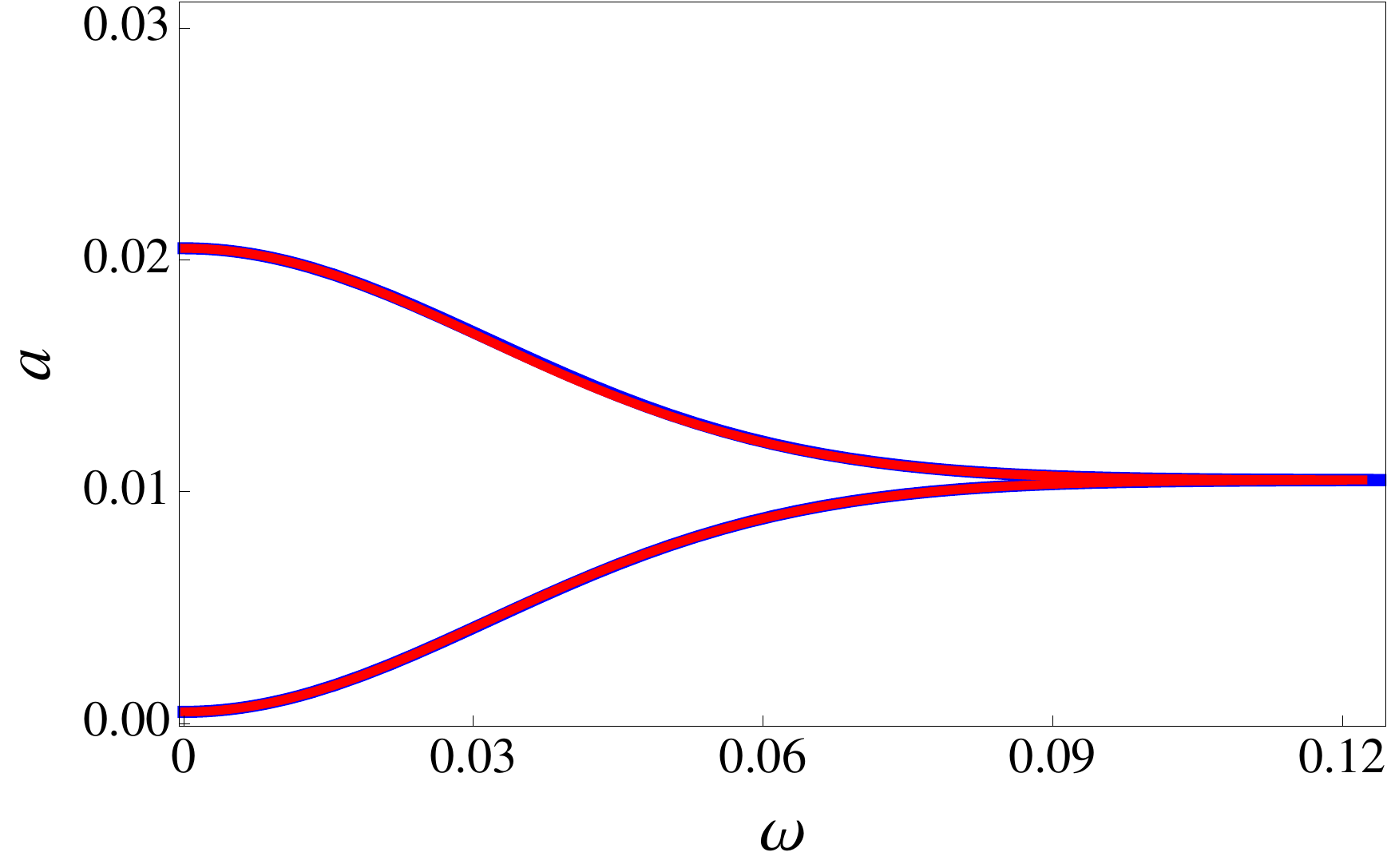}} \hspace{10pt}
\subfigure[$c=1.525$]{\includegraphics[width=2.475in]{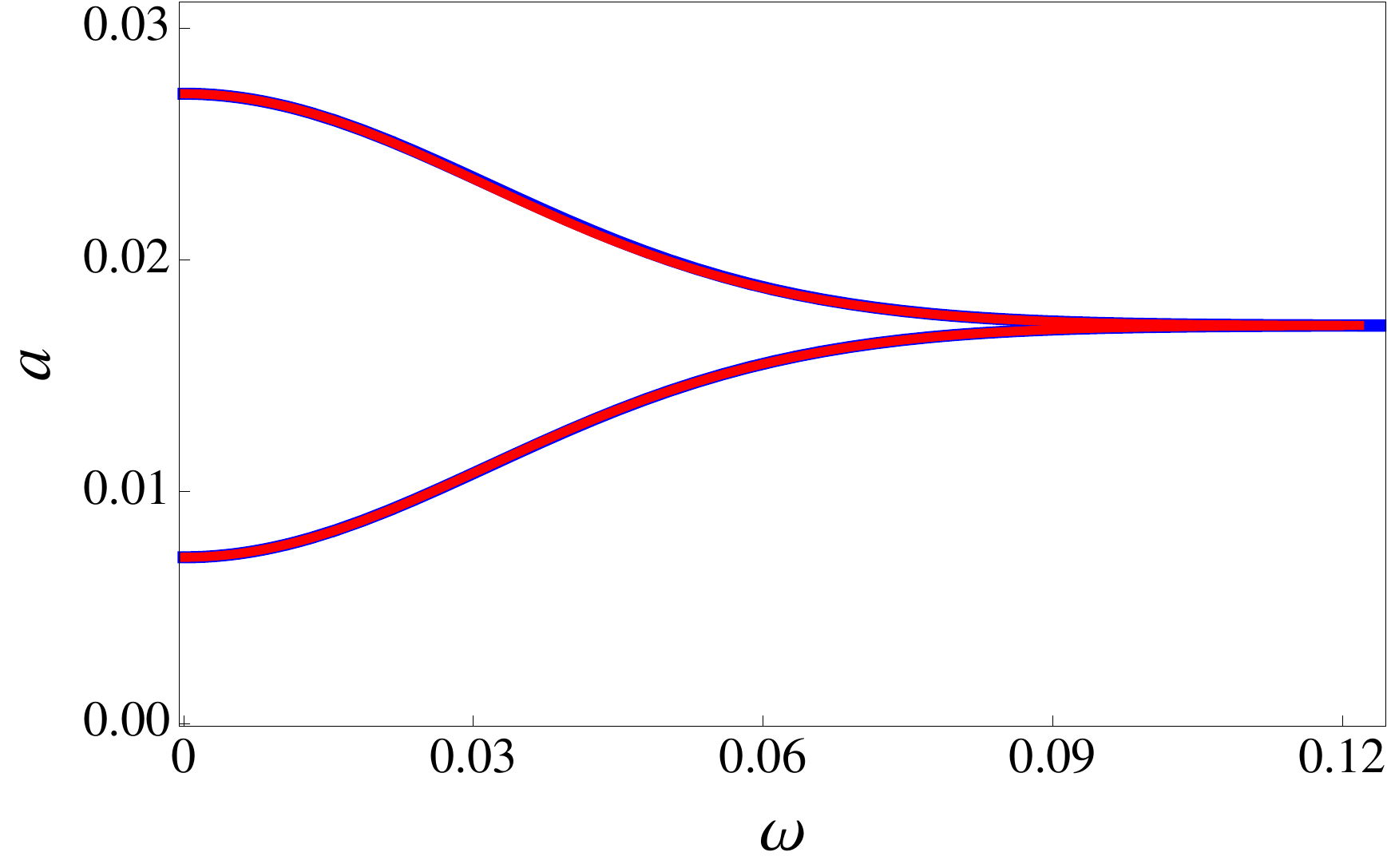}} \\
\subfigure[$b=0.015$]{\includegraphics[width=2.475in]{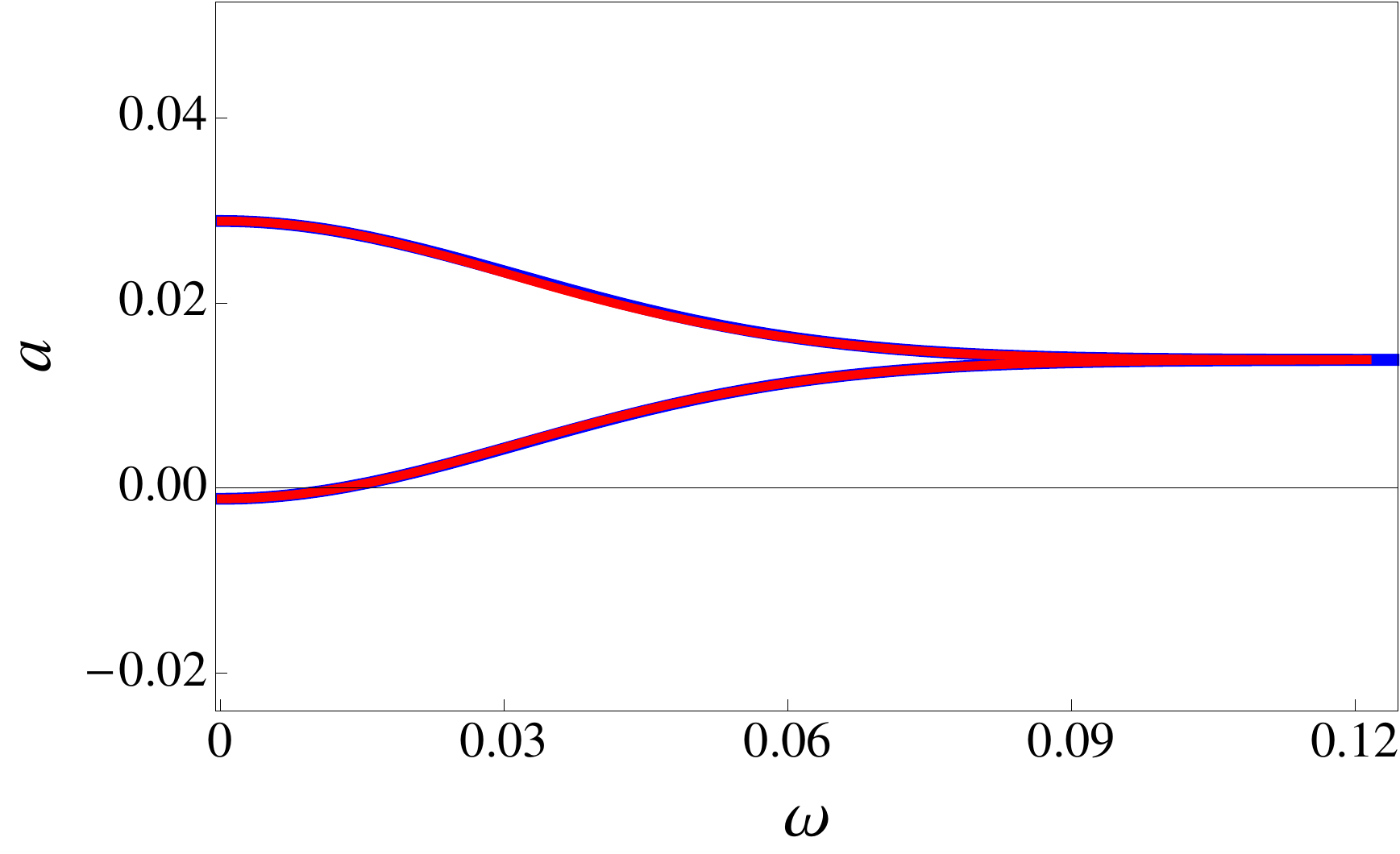}} \hspace{10pt}
\subfigure[$b=0.035$]{\includegraphics[width=2.475in]{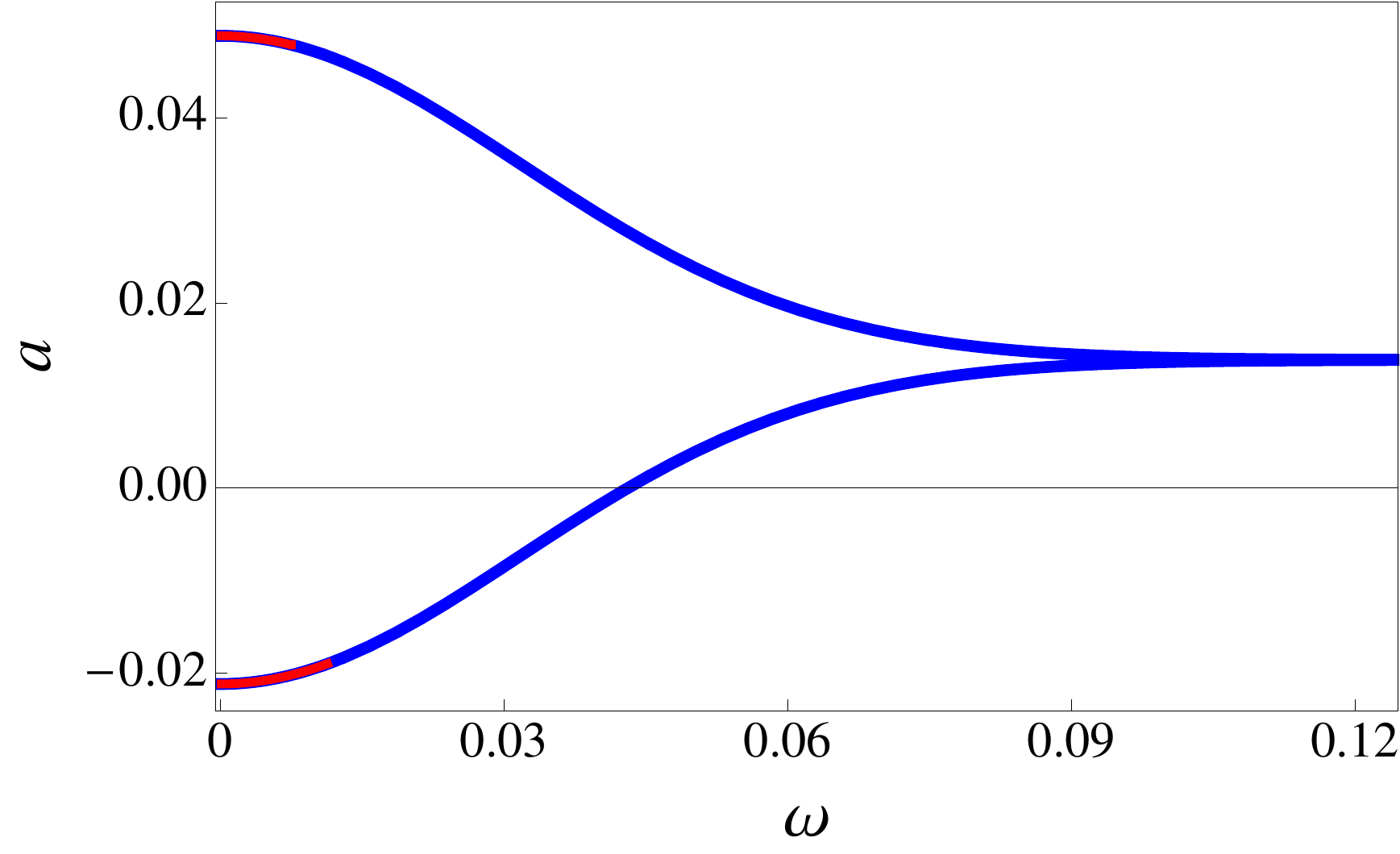}}
\caption{Comparisons between the theoretical (blue) and numerically computed (red) boundaries in the forced FitzHugh-Nagumo system \eqref{eqFH-N}. Unless stated otherwise, the parameters were fixed at $I = 0, b=0.01, c=1.52$ and $\epsilon = 0.001$. Top and middle rows: increases in $I$ and $c$ simply shift the envelope vertically up in the $(\omega,a)$ plane, in accordance with \eqref{eq:FHNenvelope}. Bottom row: the canard analysis is valid for small values of the forcing amplitude $b$, and the theory does not appear to extend for larger values. Also, the numerical continuation of the primary canards terminates (see frame (f)), and preliminary numerical continuation results reveal different dynamics there.
}
\label{fig:comparison}
\end{figure}

For low frequency forcing where $\omega = \epsilon \overline{\omega}$, we have (by Theorem \ref{thmLowFreqF}, formula \eqref{eqthmLowFF}) that primary canards of folded singularities exist in the interior of the region enclosed by
\begin{equation} 	\label{eq:FHNlow}
a = -1+ \left( \frac{2}{3}+I\right) c+\frac{\epsilon}{8} +\frac{c\epsilon}{4} \pm b \exp \left(-\frac{\epsilon \overline{\omega}^2}{2} \right). 
\end{equation}
Similarly, in the intermediate frequency forcing regime where $\omega = \sqrt{\epsilon} \Omega$, we have (by Theorem \ref{thmIntFreqF}, formula \eqref{eqthmIntFF}) that torus canards exist in the interior of the region given by
\begin{equation}  	\label{eq:FHNint}
a = -1+ \left( \frac{2}{3}+I\right) c+\frac{\epsilon}{8} +\frac{c\epsilon}{4} \pm b {\rm exp} \left(-\frac{\Omega^2}{2} \right).  
\end{equation}
We find that \eqref{eq:FHNlow} and \eqref{eq:FHNint} are different representations of the same formula
\begin{equation} 	\label{eq:FHNenvelope}
a = -1+ \left( \frac{2}{3}+I\right) c+\frac{\epsilon}{8} +\frac{c\epsilon}{4} \pm b {\rm exp}\left(-\frac{\omega^2}{2 \epsilon} \right).  
\end{equation}
This unification confirms that the folds of maximal canards of FSN I type created in the low frequency regime continue into the maximal torus canards that exist in the intermediate frequency regime. 
Figure \ref{fig:comparison} shows the region of the $(\omega,a)$ plane enclosed by the envelope \eqref{eq:FHNenvelope}.  Also, in the parameter regime of high frequency forcing with $\omega = \mathcal{O}(1)$, the results for the forced FitzHugh-Nagumo equation \eqref{eqFH-N} and for general systems \eqref{eqGen} are similar to those established in Corollary 1.3 in \cite{burke2015} for the forced van der Pol equation, and the two curves are exponentially close.



Also shown in Figure \ref{fig:comparison} are the numerically computed boundaries between the regions of parameter space which support (folded singularity and torus) canard solutions and the regions that have no canard solutions. These boundaries were computed in \textsc{auto} \cite{doedel2007} using the homotopic continuation methods developed in \cite{desroches2008,desroches2010}. More specifically, the parameters were initially chosen to lie in the interior of \eqref{eq:FHNenvelope} for low frequency forcing, and the primary maximal canard of a folded node was computed for this fixed parameter set using the methods of \cite{desroches2008,desroches2010}. We then numerically continued this canard solution in $a$ until a fold of canards (corresponding to the primary maximal canard of the FSN I) was detected. The red boundaries in Figure \ref{fig:comparison} were then obtained by switching to two-parameter continuation of this fold of canards.


For completeness, we remark that the application of Theorems \ref{thmLowFreqF} and \ref{thmIntFreqF} to the forced van der Pol system \eqref{eq:fvdp} recovers Theorems 1.1 and 1.2 from \cite{burke2015}.

\section{Li\'{e}nard Systems Subject to Small-Amplitude Forcing} \label{sec:lienard}

We conclude by showing that our results are independent of whether the forcing enters the system via a slow direction or a fast direction, provided it is small amplitude. We first verify the statement for the general class of forced Li\'{e}nard-type systems of the form
\begin{align} \label{eq:lienard} 
\frac{d^2u}{d\tau^2} + \frac{\partial f}{\partial u} (u,p,\widetilde{\epsilon}) \frac{du}{d\tau} + \widetilde{\epsilon} g\left( u, \frac{du}{d\tau},a,p,\widetilde{\epsilon} \right) = \widetilde{\epsilon}\, \widetilde{b} \cos (\omega \tau), 
\end{align}
where $f$ and $g$ are sufficiently smooth functions. Such forced Li\'{e}nard equations have two equivalent systems representations.  
In the first case, we can write the system as
\begin{equation}	\label{eq:lienardstandard}
\begin{split}
u^\prime &= v - f(u,p,\widetilde{\epsilon}), \\
v^\prime &= \widetilde{\epsilon} \left( -g\left(u,v-f(u,p,\widetilde{\epsilon}),a,p,\widetilde{\epsilon} \right) + \widetilde{b} \cos \theta \right), \\
\theta^\prime &= \omega,
\end{split}
\end{equation}
in which case the forcing term enters via a slow direction. 
%
The alternative (and equivalent) systems representation of \eqref{eq:lienard} is
\begin{equation}	\label{eq:lienardalternative}
\begin{split}
u^\prime &= v - f(u,p,\widetilde{\epsilon}) + \frac{\widetilde{\epsilon}\, \widetilde{b}}{\omega} \sin \theta, \\
v^\prime &= - \widetilde{\epsilon} g \left( u,v-f(u,p,\widetilde{\epsilon}),a,p,\widetilde{\epsilon} \right) , \\
\theta^\prime &= \omega.
\end{split}
\end{equation}
Thus, the systems representation of forced Li\'{e}nard equations of the form \eqref{eq:lienard} can have the forcing placed in either the slow or fast direction, without loss of generality. 

We observe that the Taylor series expansion of a forced planar slow/fast system near a canard point, given in \eqref{eqOLowFreqF}, is in the form \eqref{eq:lienardstandard}. That is, the dynamics of a forced slow/fast system near a canard point are described (to leading order) by the forced Li\'{e}nard equation \eqref{eq:lienard} with
\[ f(u,p,\widetilde{\epsilon}) := u^2-c_5 \widetilde{\epsilon}, \quad
g \left( u,\frac{du}{d\tau},a,p,\widetilde{\epsilon} \right) := u-c_6 u^2-c_6 \frac{du}{d\tau} + \left( c_5 c_6 - c_7 \right) \widetilde{\epsilon} - \widetilde{a}, \quad
\omega = -\frac{\widetilde{\epsilon} \, \overline{\omega}}{c_2c_3}, \]
where $c_2 c_3 <0$ (by Assumption \ref{ass:hopf}).

Finally, we note that the forced van der Pol equation \eqref{eq:fvdp} is a forced Li\'{e}nard system described by
\[ \frac{d^2x}{d\tau^2} + (x^2-1) \frac{dx}{d\tau} + \epsilon(x-a) = \epsilon b \cos (\omega \tau). \]
Similarly, the forced FitzHugh-Nagumo equation \eqref{eqFH-N} is a forced Li\'{e}nard system described by 
\[ \frac{d^2x}{d\tau^2} - (1-x^2) \frac{dx}{d\tau} + \epsilon \left( x+a-c\left( x-\frac{x^3}{3} +I - \frac{dx}{d\tau} \right) \right) = -\epsilon b \cos (\omega \tau). \]

\noindent{\bf Acknowledgments.}  The research of all three authors
was partially supported by NSF-DMS 1109587. The authors thank Mathieu
Desroches for useful conversations.

\small
\bibliographystyle{plain}
\bibliography{Draft}

\end{document}